\magnification\magstep1
\input epsf.sty
%\hsize5.2truein
%\headline{File: horseshoe10\hfill January 10, 2005}
\centerline{\bf The H\'enon Family:}
\centerline{\bf The Complex Horseshoe Locus and Real Parameter Space\footnote*{\rm Research supported in part by the NSF.}}
\medskip
\centerline{Eric Bedford and John Smillie} 

\bigskip\noindent{\bf 0. Introduction.}  The H\'enon family of diffeomorphisms 
of ${\bf R}^2$ has received a great deal of attention as a model family of dynamical systems. H\'enon was particularly interested in the question of strange attractors, but there are many other interesting questions associated with this family. Here we consider some questions concerning the global nature of parameter space.
We can write the H\'enon family of quadratic maps as $f_{a,b}:{\bf R}^2\to{\bf R}^2$, where
$$f_{a,b}(x,y)=(x^2+a-by,x).$$
For fixed parameters $a$ and $b$, $f_{a,b}$ defines a dynamical system, and we are interested in the way that the dynamics varies with the parameters. The parameter $b$ is the Jacobian of the map. When $b\ne 0$ these maps are diffeomorphisms. When $b=0$, the map has a one-dimensional image and is essentially equivalent to the one-dimensional map $g_a(x)=x^2+a$.  For $b\ne0$, Devaney and Nitecki [DN] showed that maps similar to Smale's horseshoe example occur when $a\ll0$.  On the other hand, when $a\gg0$ all points are wandering.   For $a$ in the intermediate region,  many questions concerning the  dynamics of $f_{a,b}$ are open. As an example of the type of open question we have in mind we mention the pruning front conjecture of Cvitanovi\v c. This conjecture asks whether all maps in the H\'enon family can be understood in terms of horseshoe dynamics with certain collections of orbits removed. So far evidence for this conjecture is numerical rather than theoretical. For example, work of Davis, MacKay and Sannami [DMS] can be interpreted as providing numerical evidence that for a particular value of $(a,b)$, the dynamics of the map $f_{a,b}$ can be described in terms of the horseshoe.   It is hard to imagine how one might attack this conjecture theoretically if one works solely within ${\bf R}^2$.

 The H\'enon family has a natural extension to ${\bf C^2}$. John Hubbard (see [H]) suggested in the 80's that this complex H\'enon family would also be a natural object  to study. This suggestion was undoubtedly motivated by the interesting work that was being done on the Mandelbrot set and complex dynamics in one variable.  In this note we would like to advertise some recent conjectures of Hubbard about the real and complex H\'enon families and the connection between them. The fundamental conjecture is that regions in the real parameter space which appear unrelated may in fact be have a relation which becomes apparent after we complexify both the dynamical space and the parameter space. These conjectures are intriguing because they suggest a tool which could be used to attack questions about the global structure of the real parameter space. A second purpose of this note is to show the relevance of some results from the papers [BS] and [BS{\it ii}] to Hubbard's conjectures.  We  thank Hubbard for sharing his thoughts and enthusiasm for this subject with us.

Before describing Hubbard's conjectures we will give some definitions and review in more detail some of the known results about real H\'enon diffeomorphisms. Let $\Sigma_2$ denote the space of bi-infinite sequences on two symbols, and let $\sigma$ denote the shift map on $\Sigma_2$. This dynamical system is the full 2-shift. Smale's horseshoe map of ${\bf R}^2$ has the property that the restriction to the chain recurrent set is hyperbolic and topologically conjugate to $\sigma$.  We take these as the defining properties of horseshoes and so let ${\cal H}$ be the set of (real) parameters for which the restriction of $f$ to the set of bounded orbits is hyperbolic and topologically conjugate to the two shift.

 Devaney and Nitecki [DN] show that for each $b\ne 0$ there is an $a_0(b)$ so that $(a,b)\in{\cal H}$ whenever $b<a_0(b)$.  Techniques of Hubbard and Oberste-Vorth improve on those of Devaney and Nitecki and show that we can take $a_0(b)=-2(|b|+1)^2$. When $b\ne0$ the curve $a_0$ is not in fact the boundary of the horseshoe region.  When $|b|<.08$ is small the boundary of the horseshoe locus is known to consist of curves of homoclinic ($b>0$) and heteroclinic ($b<0$) tangencies  [BS{\it ii}].  (See Figure 9.)

There are two senses in which we can complexify maps in the H\'enon family. We can complexify the dynamical space by allowing $x$ and $y$ to be complex. After complexifying the dynamical space we can also complexify the parameter space by allowing $a$ and $b$ to be complex. 
 Let us define
$$K_{a,b}=\{p\in{\bf C}^2:\{f^np:n\in{\bf Z}\}{\rm\ is\
bounded}\}.$$
Let $K_{a,b}^{\bf R}=K_{a,b}\cap{\bf R}^2$.
Thus the real horseshoe locus ${\cal H}$ be the set of (real) parameters $a,b$ for which the restriction of $f_{a,b}$ to $K_{a,b}^{\bf R}$ is hyperbolic and topologically conjugate to the full 2-shift. 

In addition to considering real horseshoes we can consider complex horseshoes. We define the complex horseshoe locus to be the set of parameters $(a,b)\in{\bf C}^2$ for which the restriction of $f_{a,b}$ to $K_{a,b}\subset{\bf C}^2$ is hyperbolic and topologically conjugate to the full 2-shift. Let us denote this by $\cal H^{\bf C}$. The results of Hubbard and Oberste-Vorth alluded to above are an outgrowth of their investigations of the complex horseshoe locus. Hubbard and Oberste-Vorth discovered that $\cal H^{\bf C}$ contains the parameter region  $HOV=\{(a,b)\in {\bf C}^2: |a|>2(|b|+1)^2, b\ne0\}$ (cf. [Ob] or [MNTU]). 
 
Real and complex horseshoes satisfy Axiom A and are structurally stable, so $\cal H^{\bf C}$ is an open set of ${\bf C}^2$, and ${\cal H}$ is an open subset of ${\bf R}^2$.  It seems likely that $\cal H^{\bf C}$ is connected, but this is not known. Let  ${\cal H}^{\bf C}_0$ be the component of $\cal H^{\bf C}$ that contains $HOV$ and let ${\cal H}_0$ be the component of $\cal H$ that contains the Devaney-Nitecki region.

The relationship between real and complex horseshoes is somewhat subtle. Consider $(a,b)$ real.  A first possibility for $(a,b)\in{\cal H}^{\bf C}$ is that $(a,b)\in{\cal H}^{\bf R}$.  A second possibility is that $(a,b)\in {\cal H}$, but $K_{a,b}\cap{\bf R}^2=\emptyset$. This occurs for example when $a>2(|b|+1)^2$. In this case, the dynamics of $f_{a,b}|{\bf R}^2$ is wandering, in fact all orbits tend to $\infty$ in both forward and backward time. Let us say that a point $(a,b)\in{\cal H}_0^{\bf C}$ is of type 3 if neither of these two possibilities occur.  That is, type 3 means: $K_{a,b}\not\subset{\bf R}^2$, and $K_{a,b}\cap{\bf R}^2\ne\emptyset$. 

\proclaim Conjecture 1 (Hubbard). There are (real) parameter values $(a,b)$ in  ${\cal H}^{\bf C}$ of type 3.

We can ask whether the particular H\'enon diffeomorphism studied by Davis, MacKay and Sannami [DMS] is a type 3 example. 
 
A real map of type 3 comes equipped with an identification of its chain recurrent set with a subshift of $\Sigma_2$.  Thus progress on Conjecture 1 might lead to insights into the Cvitanovi\v c picture of dynamics of the real H\'enon family. 

Let us fix a basepoint $(a_0,b_0)$ in ${\cal H}$ corresponding to a Devaney-Nitecki horseshoe.  There is a canonical conjugacy $h_0:K_{a_0,b_0}\to\Sigma_2$ (see \S1,2).  Let $\gamma$ be a loop in ${\cal H}_0^{\bf C}$ based at $(a_0,b_0)$.  Structural stability produces a path of conjugacies $h_t$ from the dynamical systems $f_{\gamma(t)}|K_{\gamma(t)}$ to the 1-shift.  Since $\gamma(0)=\gamma(1)=(a_0,b_0)$, the conjugacy $h_1$ induces an element $\rho(\gamma):=h_1\circ h_0^{-1}\in Aut(\Sigma_2)$.  This yields a canonical anti-homomorphism
$$\rho:\pi_1({\cal H}_0^{\bf C},(a_0,b_0))\to Aut(\Sigma_2);$$
if $\gamma_{1}\cdot\gamma_{2}$ represents the loop obtained by first following $\gamma_{1}$ and then $\gamma_{2}$, then $\rho(\gamma_{1}\cdot\gamma_{2})=\rho(\gamma_{2})\rho(\gamma_{1})$.   The image of $\rho$ is a subgroup of $Aut(\Sigma_2)$, which we denote by $\Gamma$.  Note that if we change the base point $(a_0,b_0)$ to another point $(a',b')\in{\cal H}_0^{\bf C}$, the resulting homomorphism is given by $\rho'=g^{-1}\rho g$ for some $g\in\Gamma$.  The set of parameters $HOV\subset{\cal H}^{\bf C}_0$ is not simply connected.  We will see in \S1,2 that if $\gamma$ is an element of  $\pi_1(HOV)$, then $\rho(\gamma)$ is either the identity or the automorphism of $\Sigma_2$ which interchanges the two symbols.

\bigskip
\epsfxsize=4in
\centerline{\epsfbox{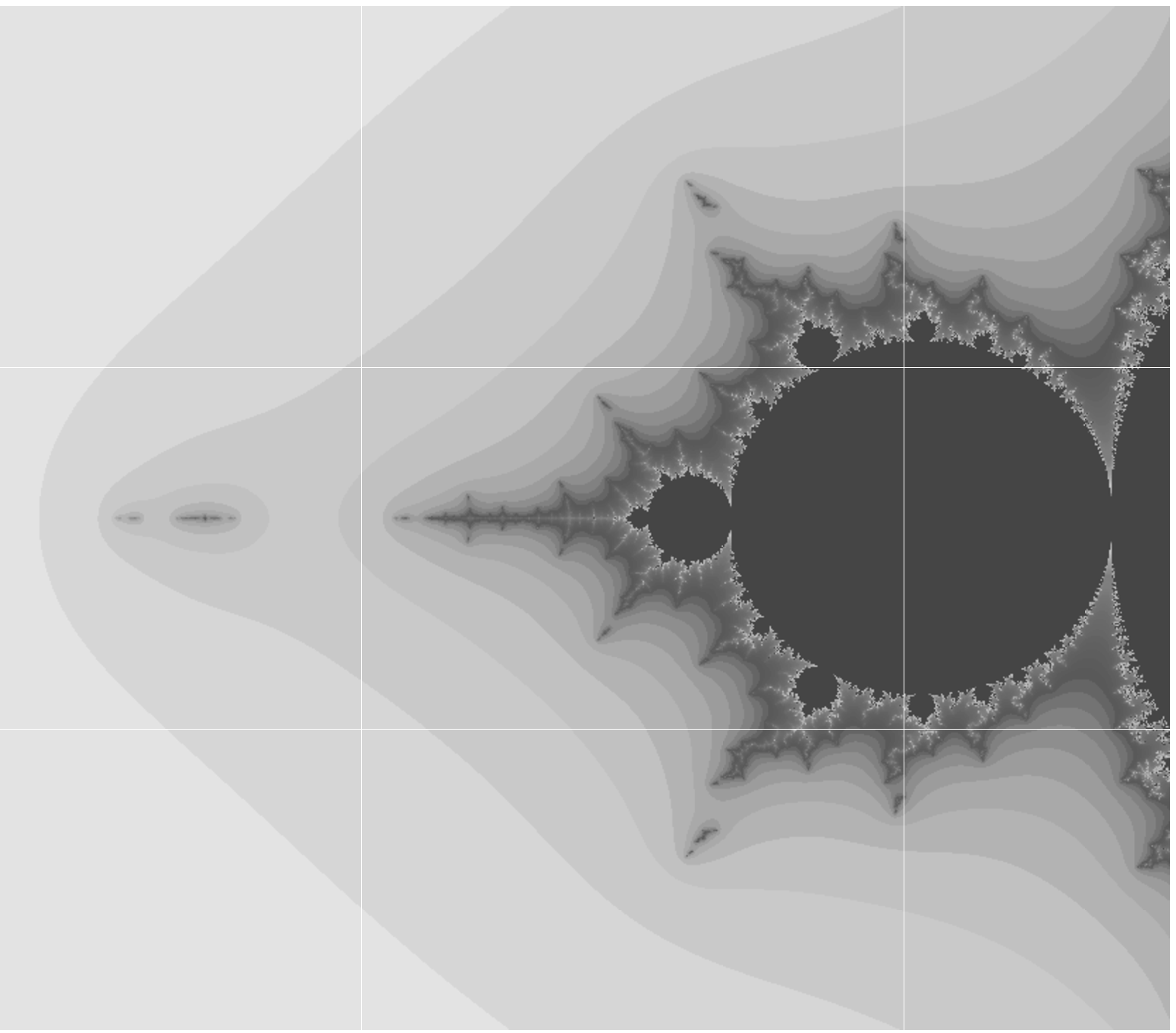}}
\centerline{Figure 1.  The $b=.2$ slice of the complex horseshoe locus.}
\medskip
The first result of this paper is to set up a correspondence between components of ${\bf R}\cap{\cal H}^{\bf C}_0$ and conjugacy classes of involutions.  From this we find:
\proclaim Theorem 1.  The number of elements of $\Gamma$ is at least as large as the number of conjugacy classes of real mappings $f_{a,b}|K^{\bf R}_{a,b}$ for $(a,b)\in {\bf R}^{2}\cap{\cal H}^{\bf C}_0$.

We have seen two conjugacy classes already: $HOV\cap\{(a,b)\in{\bf R}^2,a<0\}\subset{\cal H}$ coresponds to real horseshoes, and $HOV\cap\{(a,b)\in{\bf R}^2,a>0\}$ corresponds to the case $K^{\bf R}=\emptyset$.  The existence of a type 3 example would give a third real conjugacy type, and thus $\Gamma$ would have more than two elements.

Hubbard's conjectures are motivated in part by the computer investigations
he has carried out on the 2-shift locus in the space of quadratic
automorphisms of ${\bf C}^2$, using the computer program {\sl SaddleDrop}, which he developed with K. Papadantonakis.\footnote*{This software is available at
{\tt www.math.cornell.edu/$\sim$dynamics}; see [HP] for a detailed description.}  
Hubbard has made the intriguing
discovery of a number of ``channels'' which seem to lie inside the space of
complex horseshoes.  We use this terminology because there seem to be
``islands'' of non-horseshoe points, and it seems possible to navigate around these
islands through channels inside the horseshoe locus.  The software suggests (but does not prove) that moving around these
channels induces nontrivial automorphisms of the 2-shift.  Figure 1 gives a {\sl SaddleDrop} picture of the parameter region  $b=.2$, $-2.6<\Re(a)<-1$, $-.7<\Im(a)<.7$.  The dark region shows non-horseshoe parameters, and the lighter gray region shows points which are possible horseshoes.  Figure 2 corresponds to  $b=.2$, $-2.45<\Re(a)<-2.2$, $-.02<\Im(a)<.05$ and gives a magnification of the island on the left of Figure 1.  Many possible channels are in evidence.  The horizontal axis of symmetry in this picture is the real axis in the $b$-plane.  Points on this line that do not lie on islands are candidates for type 3 maps.
\bigskip
\epsfxsize=4in
\centerline{\epsfbox{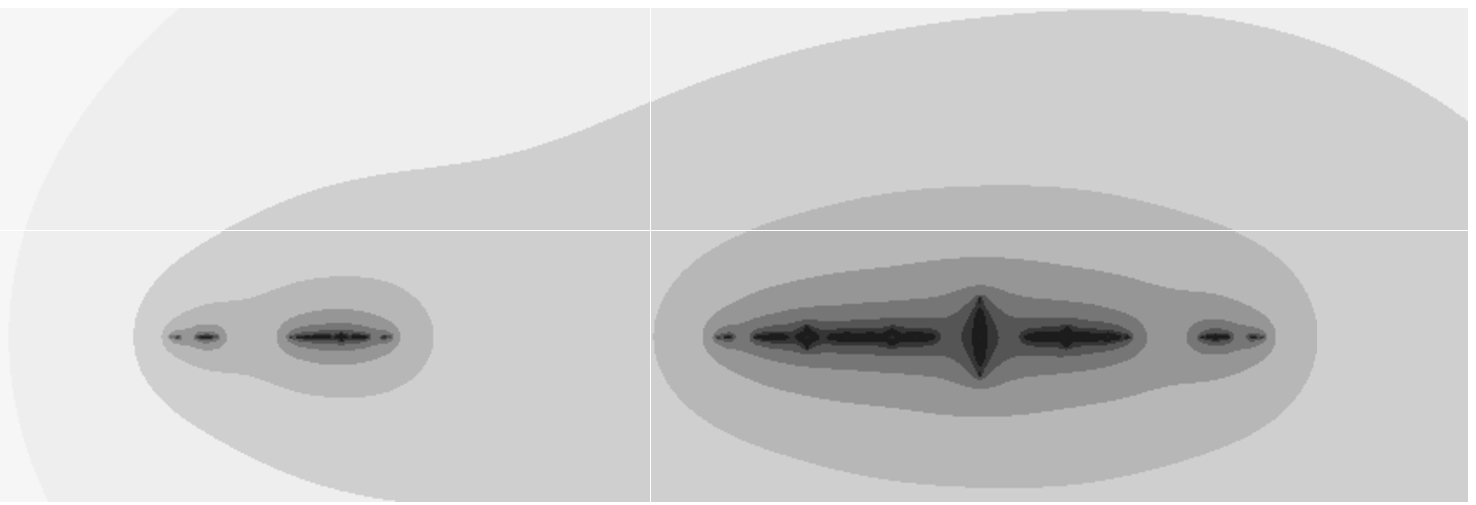}}

\centerline{Figure 2.  Detail of the ``island'' on the left side of Figure 1.}
\medskip

It is useful to compare Hubbard's observations with computer investigations which
have been carried out on real parameter space; two of these are [EM] and [HW].  In the
pictures of [EM] and [HW], there are a large number of bifurcation curves
near the boundary of horseshoe locus.   The appearance of gaps between these
bifurcation curves is consistent with the existence of places where the Hubbard's
channels  appear to intersect ${\bf R}^2$.  Such intersection points
$(a,b)\in{\cal H^{\bf C}}\cap {\bf R}^2$ would correspond to a real diffeomorphism
$f_{a,b}$ which is a complex horseshoe but (by [BS{\it ii}]) is not a real horseshoe, and thus
$K_{a,b}\not\subset{\bf R}^2$ by [BLS].

Let us remark in passing that a different sort of ``channel'' has been observed in
the work of Ricardo Oliva [Ol].  Recall that the connectivity locus is
$${\rm Conn}=\{(a,b)\in{\bf C}^2:J_{a,b}{\rm\ is\ connected}\}.$$
Clearly ${\rm Conn}\cap{\cal H}=\emptyset$.  It is also true that ${\rm
Conn}\cup{\cal H}$ is a proper subset of ${\bf C}^2$.  Oliva has found computer
evidence of ``channels'' inside the slices $\{b=const\}\cap({\bf C}^2-{\rm
Conn})$.  These channels suggest that ${\rm Conn}$ may be disconnected.

Automorphism representations similar to $\rho$ have been considered for polynomial maps of ${\bf C}$ of degree two and three. The analog of the horseshoe locus is the shift locus. In the degree 2 case the celebrated result on the connectivity of Mandelbrot set shows that the shift locus has infinite cyclic fundamental group. The automorphism group of the one sided 2 shift is cyclic of order two and the homomorphism is surjective. In degree 3 the situation is much more complicated.  Nevertheless the results of   [BDK1] (see also [BDK2]) show in this case also the corresponding map is surjective. These results can be viewed as motivation for the following conjecture.

\proclaim Conjecture 2 (Hubbard). $\Gamma=Aut(\Sigma_2)$.

Showing that $\Gamma$ is large would have
interesting consequences for the family of H\'enon diffeomorphisms. The
automorphism group of the two-sided $d$-shift is quite complicated for all values
of $d$ unlike the automorphism group of the one sided $d$-shift. In this case the
family of H\'enon diffeomorphisms (for which the degree is 2) would exhibit some
behavior quite different from that exhibited by the family of quadratic maps of ${\bf C}$.

In this paper we will discuss the parameter region ${\cal W}$, which is defined in [BS].  For $(a,b)\in{\cal W}$ the points of $K$ may be given symbolic codings:  this coding is based on the graph ${\cal G}$ in Figure 4, whereas the horseshoe coding is based on the graph ${\cal E}$.  This coding does not determine the dynamics, but it gives us a language with which we may identify locations within $K$ on which certain dynamical behaviors may (or may not) occur.  This language is used in [BS] to give a criterion for when $f|K^{\bf R}$ is a real horseshoe.
${\cal W}$ is a fairly substantial region of parameter space: the slice ${\cal W}\cap\{b=0\}$ is pictured in Figure 6, and ${\cal W}\cap{\bf R}^2$ is pictured in Figure 9. 

\proclaim Theorem 2.  The subgroups $\rho(\pi_1(HOV))$ and
$\rho(\pi_1({\cal W}\cap{\cal H}^{\bf C}))$, together with the shift, generate
a proper subgroup of $Aut(\Sigma_2)$.

The parameter region pictured in Figure 2 is contained in ${\cal W}$.  By Theorem 2, automorphisms corresponding to loops in Figure  2 do not generate $Aut(\Sigma_2)$.  Since ${\cal W}$ does not contain all of ${\cal H}^{\bf C}$, Theorem 2 does not resolve Conjecture 2 one way or the other, but it does imply that in order to prove Conjecture 2, it is necessary to consider paths that are not contained in ${\cal W}$.

\bigskip

\noindent{\bf 1.  One-Dimensional Horseshoes.}  When $b=0$,
the map $f_{a,b}= f_{a,0}(x,y)=(x^2+a,x)$ is not a
diffeomorphism.  The dynamical behavior of this map is determined
by the behavior of the 1-dimensional map $g_a(x)=x^2+a$.  
We begin by focussing on this one-dimensional case.  We say that $g_a$ is a complex horseshoe if it is hyperbolic (which is equivalent to
being uniformly expanding on its Julia set $J$) and conjugate to
the 1-sided shift on two symbols.  

Let ${\cal P}_1=\{(a,0):a\in{\bf C}\}$ be the set of parameters for which $f$ is not a diffeomorphism, and let $ {\cal P}_2=\{(a,b)\in{\bf
C}^2:b\ne0\}$.  
Thus ${\cal P}_1\cup{\cal P}_2={\bf C}^2$.  The complex
(one-dimensional) horseshoe locus in ${\cal P}_1$ is the
complement of the Mandelbrot set, i.e., $({\bf C}-M)\times\{0\}$.   The complex horseshoe parameters which are real are ${\bf R}\cap ({\bf C}-M)=(-\infty,-2)\cup(1/4,\infty)$.  For $a\in(-\infty,-2)$, we have $J_a\subset{\bf R}$, and thus the horseshoe is real.  On the other hand, for $a\in(.25,\infty)$ we
have $J_a\cap{\bf R}=\emptyset$.  

Let us discuss the relationship between  ${\cal P}_1$ and ${\cal H}^{\bf C}$.  For $a\in{\bf C}$, we may consider the natural extension $\hat g_a:\hat J_a\to \hat J_a$.  Here $\hat J_a$ refers to the bi-infinite sequences $\{z_n:g_a(z_n)=z_{n+1}, z_n\in J_a,n\in{\bf Z}\}$, and $\hat g_a$ acts as the shift on $\hat J_a$.  It was shown in [HO2] that if $g_a$ is hyperbolic, then there is a number $r=r(a)$ such that if $b\ne0$ satisfies $|b|<r(a)$, $f_{a,b}|J_{a,b}$ is conjugate to the natural extension $\hat g_a|\hat J_a$.  
We make two sorts of conclusions from this.  First, if $\gamma:\theta\mapsto(a,re^{2\pi i\theta})$ denotes the curve in ${\cal P}_2$ which ``winds around'' a point $a\in {\bf C}-M$, then $\rho(\gamma)$ is the identity.  (This is because the conjugacy between $f_{\gamma(\theta)}|J_{\gamma(\theta)}$ and  $\hat g_a,\hat J_a$ varies continuously in $\theta$.)   On the other hand, if $a\in M$ is a hyperbolic parameter, we conclude that there is a neighborhood $U$ of $(a,0)$ such that $U\cap {\cal H}^{\bf C}=\emptyset$.   Thus if $(a,b_0)\in{\bf R}^2$ is a point of type 3, and if $b_0$ is close to 0, then $a$ must be close to a non-hyperbolic point of $M\cap {\bf R}$.

\bigskip
\epsfysize=1.6in
\centerline{\epsfbox{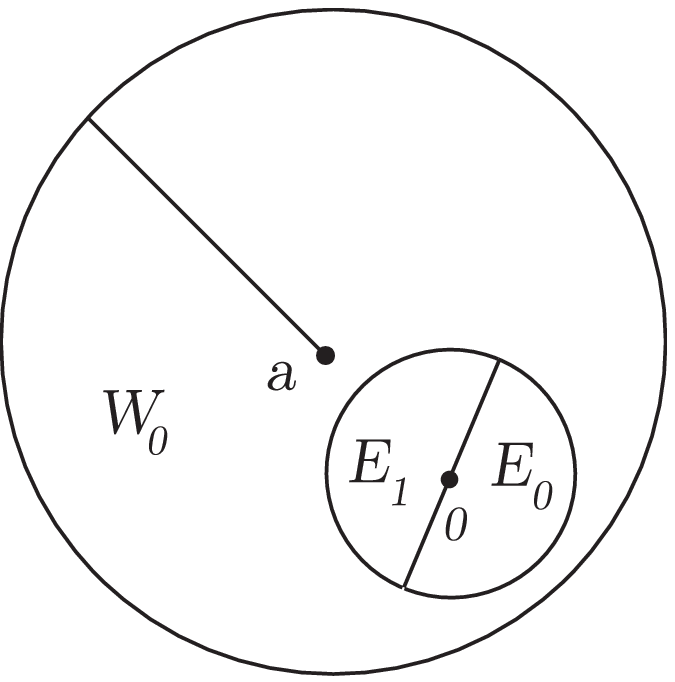}}
\centerline{Figure 3.  Construction of $E_0$ and $E_1$.}
\medskip

Abusing notation, we let $HOV=\{(a,0):|a|>2\}\subset{\cal P}_1$.
Let us recall the proof that $g_a$ is a complex
horseshoe for $a\in HOV$. The map $g_a$ takes the
disk $|z|<R$ to the disk with center $a$ and radius $R^2$.   Let $W_0$ denote the disk $\{|z-a|<R^2\}$ with a radial slit removed as in Figure 3.   Let us choose  $0\le\phi\le 2\pi$ such that $a=|a|e^{i\phi}$ and define the two half disks:
$$E_0=\{0<|z|<R:\phi/2+\pi <arg(z)<\phi/2+2\pi\}$$
and
$$E_1=\{0<|z|<R:\phi/2<arg(z)<\phi/2+\pi\}.$$
If $2<R<|a|$, then $\bar E_0\cup\bar E_1\subset W_0$.    In this case there is a coding map $c:J_a\to \Sigma_2^+$, where
$c(p)$ is the sequence whose
$n$th element is defined to be the value $k$ for which 
$g^n(p)\in E_k$.  (Here we let $\Sigma_2^+$ denote the one-sided shift space on two symbols.)  This coding map is a conjugacy.
\bigskip
\epsfxsize=4in
\centerline{\epsfbox{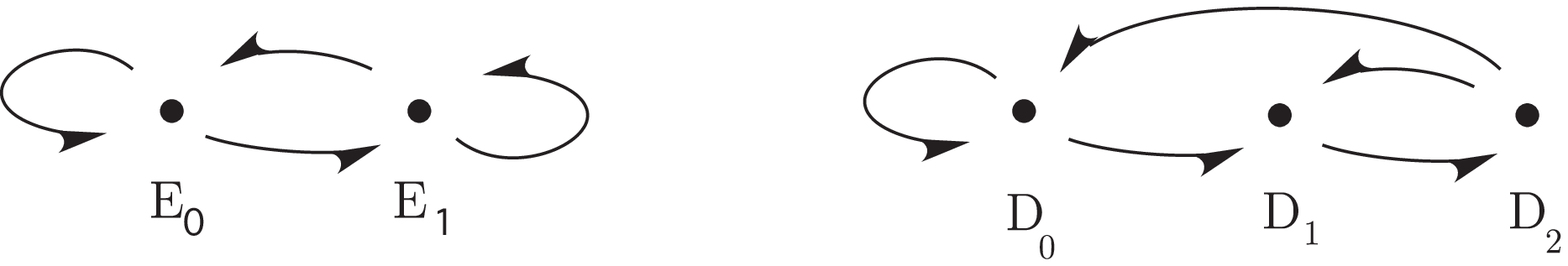}}
\centerline{Figure 4.  Graphs ${\cal E}$ (on left) and ${\cal G}$ (on right).}
\medskip

The coding that we have constructed depends continuously
on the parameter $a$. If $a$ is positive and real, we may write $a=|a|e^{2\pi i\phi}$ for $\phi=0$ and $\phi=2\pi$.  This gives two different choices of labels for the partition $\{E_0,E_1\}$.  This implies that in the definition of $c$ above, we obtain the two conjugacies to the shift which differ by an interchange of the symbols $0$ and $1$.  If we move
$a$ through a closed loop $\tau$ which winds once around the origin, then
the induced representation $\rho(\tau)\in Aut(\Sigma^+_{\cal E})$
corresponds to interchanging the symbols $0$ and $1$.

The point $a=-2$ is of particular importance since it is the boundary of
the set of real, 1-D horseshoe parameters.  We will
describe a different scheme for coding points which works for $a$ in a
region containing $-2$.  This scheme utilizes the 3-box construction
given in [BS{\it ii}].    We start by setting $D_0=E_0$, and $D_2=E_1$, with $E_j$ defined as above.
Now we choose a constant $c\in te^{i(\phi-\pi)/2}$, $t>0$, and let $W_1$ to be the portion of the disk $\{|z-a|<R^2\}$ lying to the left of $c+{\bf R}e^{i\phi/2}$.  In case $\phi=\pi$, we have $c>0$, which is pictured in Figure 5.  We define $D_1:=g_a^{-1}W_1$ to be the preimage of $W_1$.  It follows that we have proper, unbranched maps
$$D_0\to W_0,\ \ D_2\to W_0, \ \ D_1\to W_1.$$
The condition for $a\in{\cal W}$ is that
$$\bar D_0\cup\bar D_1\subset W_0,\ \ \bar D_2\subset W_1.$$
The admissible transitions $g_a:D_i\to D_j$ are given by the graph ${\cal G}$ in Figure 4.   The transition $D_i\to D_j$ being admissible means that $\bar D_j\subset W_k$, and $g_a:D_i\to W_k$ is a proper map.  The admissible
maps all have degree one, except the map $D_1\to D_2$, which has
degree 2.  Let $\Sigma_{\cal G}^+$ denote the space of 1-sided
sequences on the symbols $\{0,1,2\}$ which are consistent
with ${\cal G}$.  

\bigskip
\epsfysize=1.4in
\centerline{\epsfbox{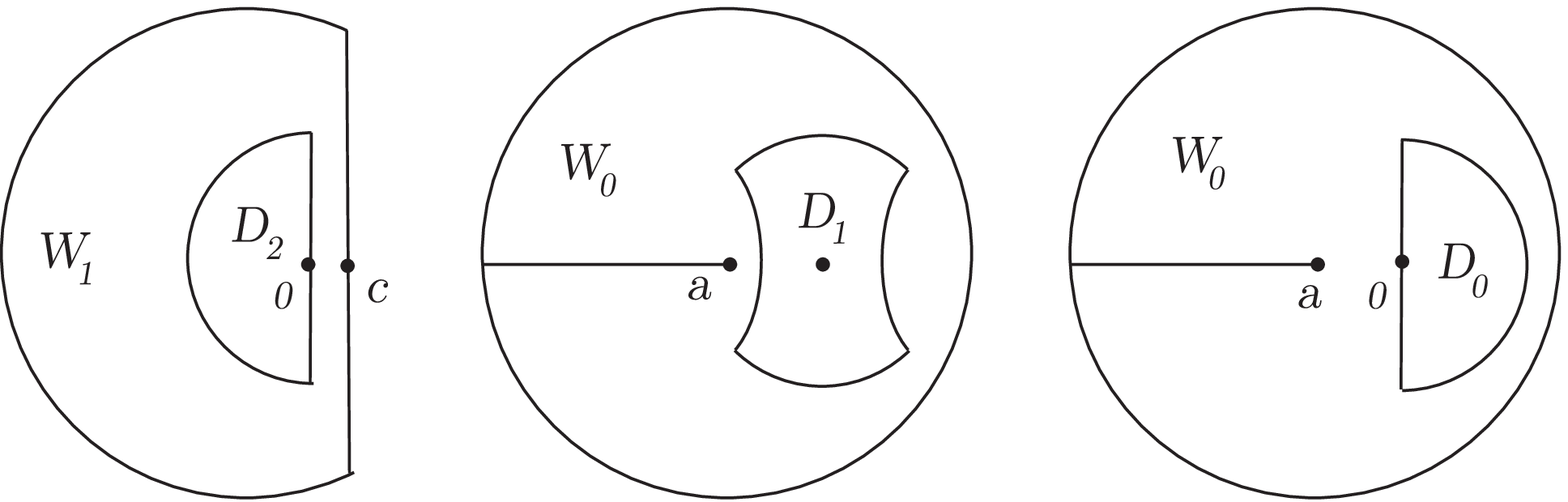}}
\smallskip
\centerline{Figure 5.  Construction of $W_j$, $j=0,1$ and  $D_k$, $k=0,1,2$.}
\medskip

In order to describe a domain ${\cal W}$ for which this construction can be carried out, let us simply consider the case $c=0$.  For a parameter $a=re^{i\phi}\in {\bf C}$, let $C_a$ denote the preimage under $g_a$ of the line ${\bf R}e^{i\phi/2}$.  The
preimage $C_a$ consists of two arches of a hyperbola.  The limiting case of the construction of the cover $\{D_0,D_1,D_2\}$ is when $a\in C_a$.  An easy calculation shows that the condition $a\in C_a$
corresponds to $a^2+a\in \sqrt a{\bf R}$, which defines an arch
of a hyperbola passing through $-1$ and with asymptotes of angle
$2\pi/3$ and $4\pi/3$.  So let us define ${\cal
B}^3$ as the connected parameter region bounded by this
hyperbola, and which contains the interval $(-\infty,-1)$.   Figure 6 shows the regions $HOV$ and ${\cal W}$, with the Mandelbrot set drawn in for comparison.  $HOV$ is the region outside the circle of radius 2, and ${\cal W}$ is the region to the left of the arch of the hyperbola.

For $(a,0)$ in ${\cal W}$, the Julia set $J_a$ is contained in $D_0\cup D_1\cup D_2$.  
For $p\in J_a$, we would like to code the forward orbit $\{g^np:n\ge0\}$, but the definition of its ``itinerary'' is not immediately clear since we could have $g^np\in D_i\cap
D_j$ for $i\ne j$.  However, by [BS{\it ii}], there is a set of codings $c_{\cal G}(p)\subset
\Sigma_{\cal G}$ with the property that if $s=s_0s_1\dots\in c(p)_{\cal G}$, then $g^np\in
D_{s_n}$ for every $n\ge0$.  The $\alpha$ fixed point $p_\alpha$ belongs to $J_a\cap
D_1\cap D_2$; both $(12)^\infty$ and $(21)^\infty$ are codings for the
forward itinerary that are consistent with the graph ${\cal G}$.  These are the only
possible such codings for $p_\alpha$.  Furthermore,
$c_{\cal G}(p)$ is a single point unless $p\in g^{-n}p_\alpha$ for some
$n\ge0$, in which case $c_{\cal G}(p)$ contains exactly two sequences, one
ending in
$(12)^\infty$ and the other ending in $(21)^\infty$.  

Note that this coding construction works in a parameter
region which contains non-hyperbolic
mappings and bifurcations; while we will obtain symbolic codings for
our mappings, we cannot expect these codings to be
conjugacies.
\bigskip
\epsfysize=2.2in
\centerline{\epsfbox{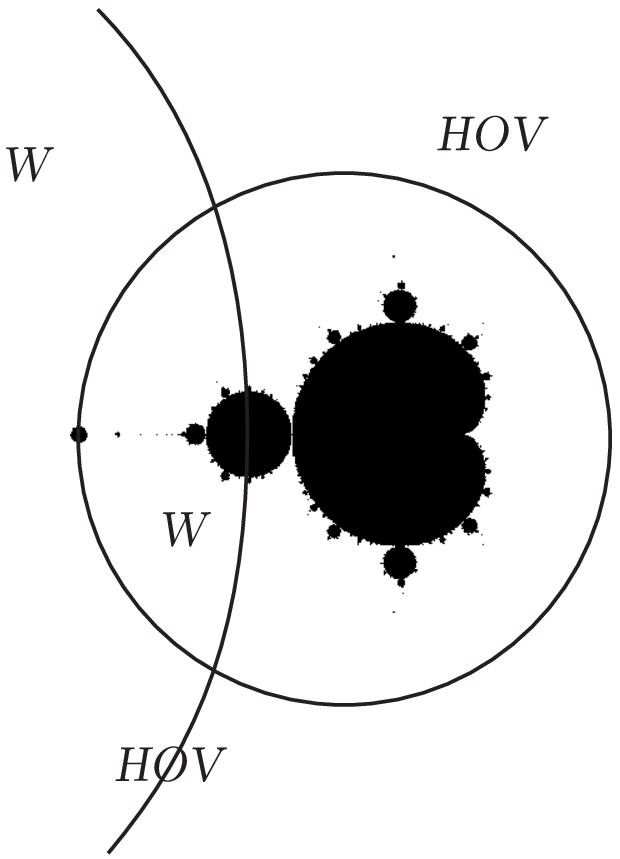}}
\centerline{Figure 6.  $HOV$, ${\cal W}$, and $M$}
\medskip

We can also describe the ${\cal W}$ codings in terms of certain sequences of
0's and 1's.  This will be the itinerary of a point with respect to the cover $\{W_0,W_1\}$.  The possible transitions for this cover correspond to the graph ${\cal G}_0$ on the left hand side of Figure 7.  This gives a coding map $c_{{\cal G}_0}:J\to\Sigma^+_{{\cal G}_0}$, where $\Sigma^+_{{\cal G}_0}$ denotes the one-sided shift space on the graph ${\cal G}_0$.  Specifically, in Figure 5, we see that $D_0$ is contained in $W_0$, and $g_a$ maps $D_0$ to $W_0$.  Thus $D_0=W_0\cap g_a^{-1}W_0$, so we may denote $D_0$ by the symbol {\bf 00}.  Similarly, $D_1$ is contained in $W_0$ and maps to $W_1$; we denote $D_1$ by the symbol {\bf 01}.  Finally, $D_2=W_1\cap g^{-1}_aW_0$, corresponding to {\bf 10}.  The re-coding of ${\cal G}$ is given on the right hand side of Figure~7.  Note that this is also the 2-block extension of ${\cal G}_0$.
\bigskip
\epsfxsize=4in
\centerline{\epsfbox{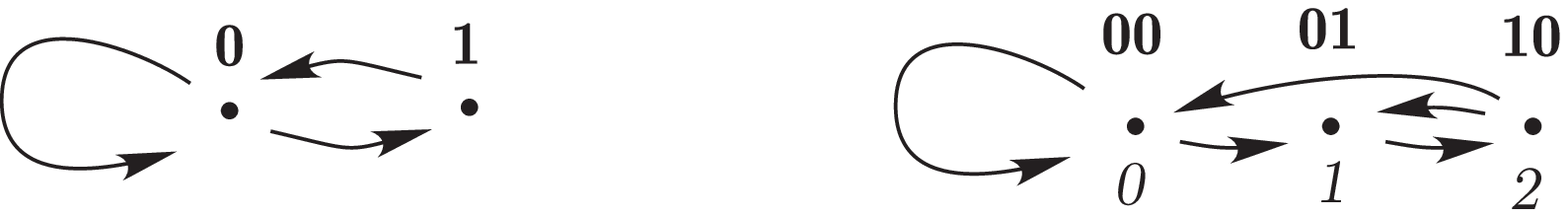}}
\smallskip
\centerline{Figure 7.  The graph ${\cal G}_0$ (on left);  re-coding on 2 symbols.}
\medskip
Finally, let us consider a parameter $a$ for which $(a,0)\in HOV\cap{\cal W}$.  In this case, we have two ways of coding points: using $\Sigma^+_{\cal E}$ and $\Sigma^+_{\cal G}$.  In order to compare these two codings, we consider the new cover obtained as the common refinement of ${\cal E}$ and ${\cal G}$:  $\{E_i\cap D_j:i=0,1,\ j=0,1,2\}$.  We have $E_0\cap D_2=\emptyset$ and $E_1\cap D_0=\emptyset$.  We denote $E_i\cap D_j$ by the symbol ${{\bf i}\atop j}$, and the graph corresponding to the admissible transitions is given in Figure~8.  It follows that the $g_a$ itinerary of a point can be assigned a coding in  $\Sigma^+_{\hat{\cal G}}$, the set of one-sided sequences consistent with ${\hat{\cal G}}$.    We have two projections:  $\pi_{\cal E}:\Sigma^+_{\hat{\cal G}}  \to\Sigma^+_{\cal E}$ given by taking the sequence of symbols on the top, and $\pi_{\cal G}:\Sigma^+_{\hat{\cal G}}\to\Sigma^+_{\cal G}$, given by taking the bottom symbols.  The map $c_{\cal E}:J\to\Sigma_{\cal E}^+$ is an (invertible) homeomorphism. The passage from an ${\cal E}$-coding to a ${\cal G}$-coding of a point $p\in J_a$ is obtained by passing through the graph $\hat{\cal G}$, following $\pi_{\cal G}\circ \pi_{\cal E}^{-1}$.
\bigskip
\epsfxsize=4in
\centerline{\epsfbox{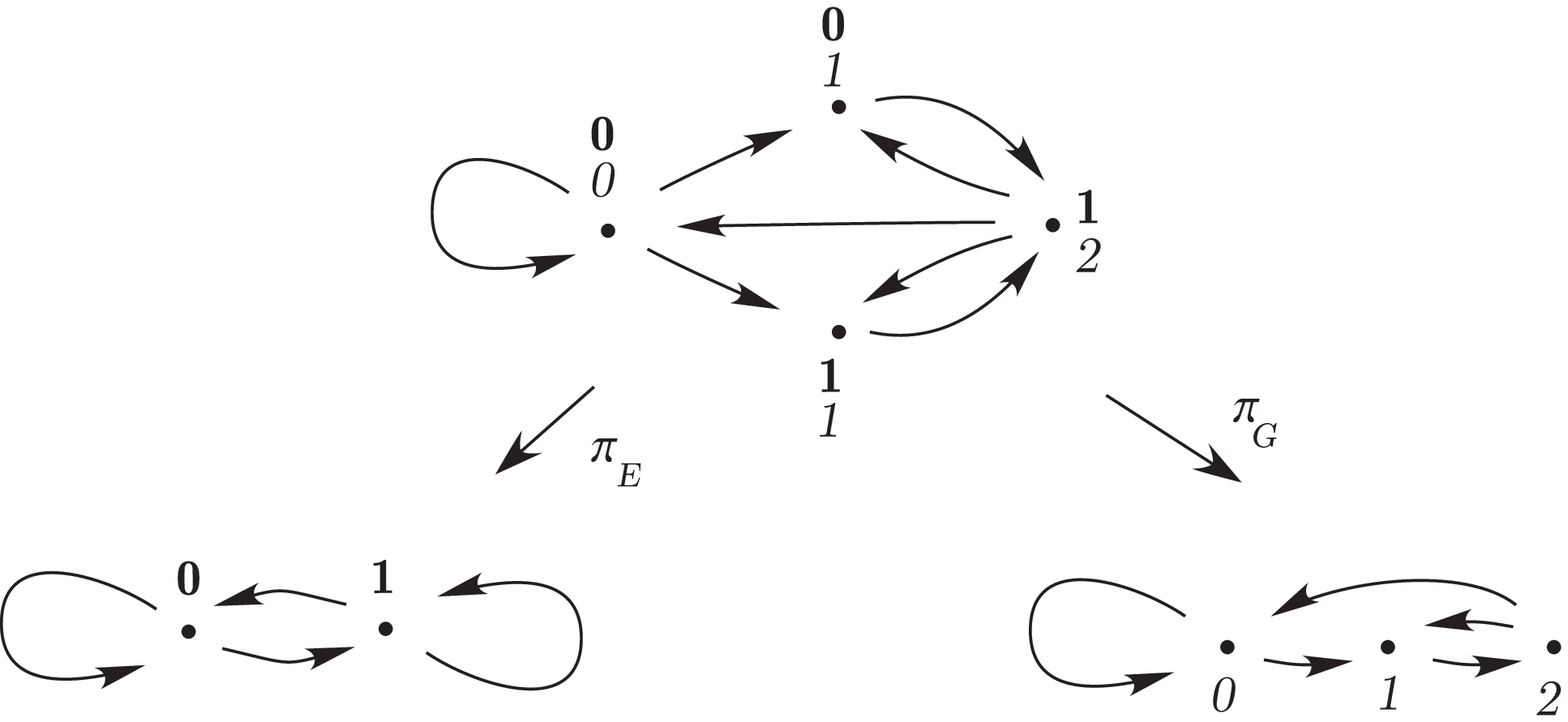}}
\smallskip
\centerline{Figure 8.  The graph $\hat{\cal G}$;  passage between codings.}
\medskip

\bigskip\noindent{\bf 2. The 2-Dimensional Case.}  For $(a,b)\in{\cal P}^2$ and $(x,y)\in{\bf C}^2$ let us write $(x_n,y_n)=f^n(x,y)$.  Then each $(x,y)\in{\bf C}^2$ gives rise to a unique bi-infinite sequence $(y_n)_{n\in{\bf Z}}$, which satisfies $y_{-1}=x$, $y_0=y$, and 
$$y_{n+1}=g_a(y_n)-by_{n-1}.$$
Recall that a sequence $(z_n)_{n\in{\bf Z}}$ is said to be an $\epsilon$ pseudo-orbit for $g_a$ if ${\rm dist}(z_{n+1},g_az_n)<\epsilon$ for all $n\in{\bf Z}$.  For each $(a,b)\in{\cal P}^2$ there is an $R=R(a,b)$ such that $K_{a,b}\subset\{|x|,|y|\le R\}$.  It follows that $(y_n)_{n\in{\bf Z}}$ is an $\epsilon$ pseudo-orbit if $\epsilon>R|b|$.  Abusing notation slightly, we define ${\cal W}$ to be the set $\{(a,b):a\in{\cal W}, |b|<\epsilon R^{-1}\}$, with $\epsilon$ and $R$ as above.  (Context should make it clear whether we are talking about the 1-dimensional or 2-dimensional versions of ${\cal W}$.)   In [BS] we estimate the value of $\epsilon$ for which the construction works, and we show that the real slice of ${\cal W}$ can be taken as the region to the left of the boundary curve which is indicated in Figure~9.  For comparison purposes, Figure 9 also gives the $HOV$ locus, as well as two darker, dashed curves, which represent portions of two curves of homoclinic/heteroclinic tangencies.  The actual curves of tangencies, in fact, continue smoothly through the point $(-2,0)$; the portions we have drawn here are the conjectural boundaries of the horseshoe locus ${\cal H}$.  In fact, this conjecture has been verified for $|b|<.08$ (see [BS{\it ii}]):  ${\cal H}\cap\{|b|<.08\}$ coincides with the region to the left of these curves. 

For $a\in{\cal W}$, let $\{W_0,W_1\}$, $\{D_0,D_1,D_2\}$, and ${\cal G}$ be as in the previous section.   Let ${\rm Orb}(\epsilon)$ denote the $\epsilon$ pseudo-orbits for $g_a$ which are contained in $D_0\cup D_1\cup D_2$.  We also consider pseudo-orbits that are equipped with codings;  let 
$$\hat X(\epsilon)=\{(p_n,i_n)_{n\in{\bf Z}}:(p_n)\in{\rm
Orb}(\epsilon), (i_n)\in\Sigma_{\cal G}, p_n\in{\cal D}_{i_n}\}.$$
We have two natural projections: $\pi_1:\hat X(\epsilon)\to{\rm
Orb}(\epsilon) $ and $\pi_2:\hat X(\epsilon)\to\Sigma_{\cal G}$.
The following result on the existence and  stability of codings for pseudo-orbits is proved in [BS]:
\proclaim Theorem 3.  For $\epsilon>0$ sufficiently small,  $\pi_1$ is a
surjection.  Further, there is a $\delta>0$ such that if $(p_n,i_n)_{n\in{\bf
Z}}\in \hat X(\epsilon)$, and if $(q_n)\in{\rm Orb}(\epsilon)$ is
$\delta$-close to $(p_n)$, then $(q_n,i_n)_{n\in{\bf Z}}\in  \hat
X(\epsilon)$.
\bigskip
\epsfxsize=2.6in
\centerline{\epsfbox{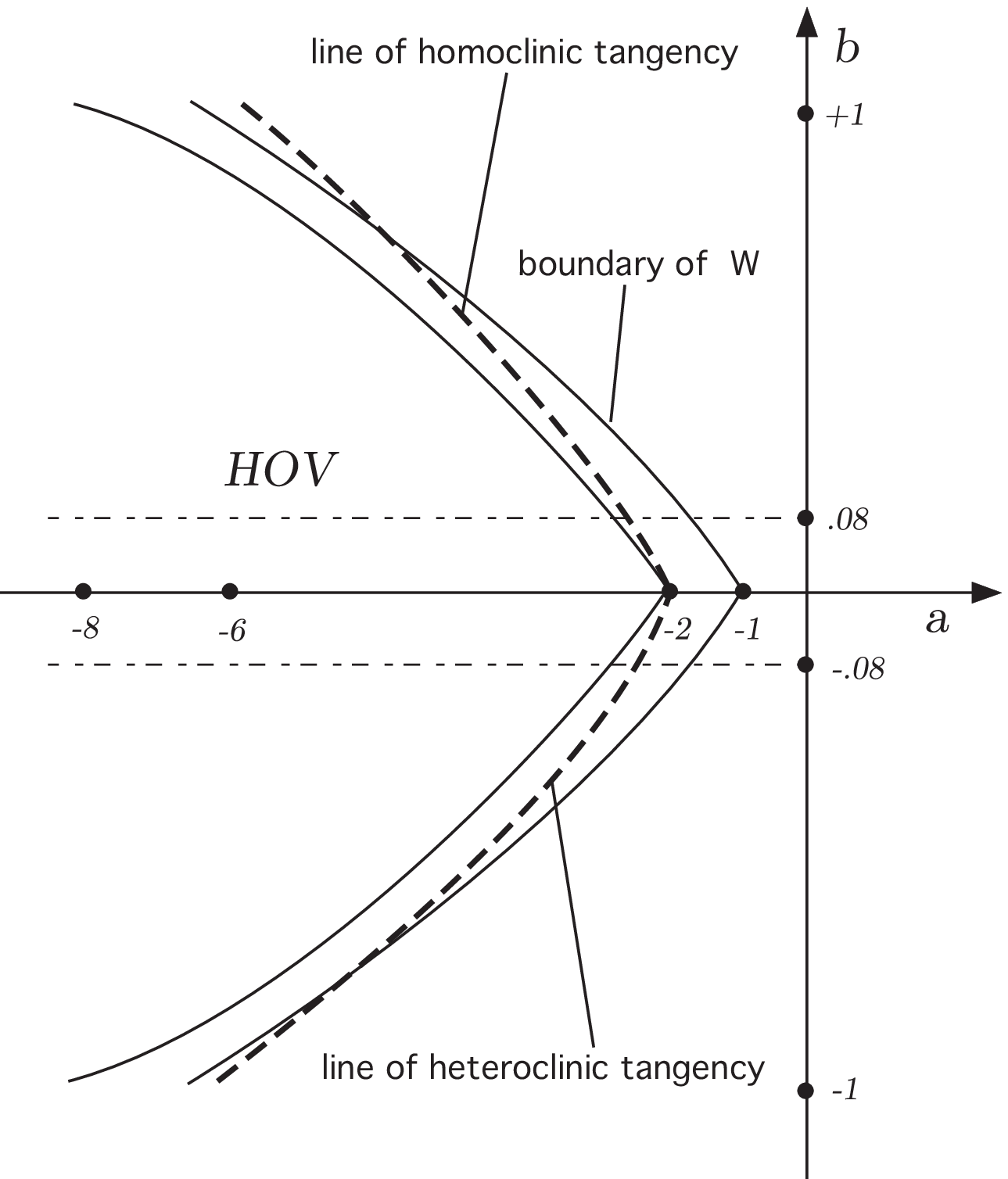}}
\smallskip
\centerline{Figure 9.  The real slice ${\bf R}^2\cap{\cal W}$}
\bigskip

Thus an $\epsilon$ pseudo-orbit $(p_n)$ corresponds to a set of symbol sequences  $c_{\cal G}:=\pi_2\circ\pi_1^{-1}(p_n)\subset  \Sigma_{\cal G}$.  We have seen that if $|b|<\epsilon R^{-1}$, then the $f_{a,b}$ orbit of a point $p\in K_{a,b}$ may be identified with an $\epsilon$ pseudo-orbit for $g_a$.  Thus we may assign these codings $c_{\cal G}(p)$ to any point $p \in K_{a,b}$.

Next we discuss how many points are in the set
$c_{\cal G}(p)$ for given $p\in J$.  One property is that for $(a,b)\in{\cal W}$,
$f_{a,b}$ has two distinct saddle points $P_\alpha\in {\cal D}_1\cap {\cal D}_2$ and
$P_\beta\in{\cal D}_0$.  These correspond in a natural way to the $\alpha$ and
$\beta$ fixed points of
$g_a$ for $(a,0)\in{\cal W}$.  
\proclaim Theorem 4.  A point $p\in J$
will have a unique coding $c_{\cal G}(p)\in\Sigma_{\cal G}$ unless $p\in
W^s(P_\alpha)$, in which case $c_{\cal G}(p)$ consists of two elements of
$\Sigma_{\cal G}$, one ending in $(12)^\infty$ and one ending in $(21)^\infty$.

Each $\gamma\in\pi_1({\cal W}\cap{\cal H^{\bf C}})$ yields a monodromy
homeomorphism
$\chi_\gamma$ of
$J_{a_0,b_0}$.  Theorem 5 gives us restrictions on the homeomorphisms
$\chi_\gamma$ that can arise.

\proclaim Theorem 5.  If $\gamma$ is a closed loop in ${\cal W}\cap{\cal
H}^{\bf C}$, then
$c_{\cal G}=c_{\cal G}\circ\chi_\gamma$.

\noindent{\it Proof.}  Let
$s\mapsto\gamma^{(s)}$ be a homotopy inside $\tilde{\cal W}$ taking
$\gamma$ to the base point
$(a_0,b_0)$.  For $(a,b)\in\gamma^{(s)}$, let 
$${\rm Per}_N(a,b):=\{p\in{\bf
C}^2:f^Np=p\}$$
denote the set of periodic points whose periods divide $N$.  Let $S:=c_{\cal
G}({\rm Per}_N(a_0,b_0))$ denote the set of admissible codes corresponding to these
periodic points.  There are only finitely many such periodic points, and each point
can have at most two distinct codings, so $S$ is a finite set.  The set
${\rm Per}_N(a,b)$ varies continuously with respect to $(a,b)$. If $(a',b')$ is
close to
$(a,b)$, then for $p'\in{\rm Per}_N(a'b')$, then $\{f^j_{a'b'}(p'):j\in{\bf
Z}\}\in{\rm Orb}_{a,b}(\epsilon)$.  Further, the $\epsilon$-orbit
$\{f^j_{a',b'}p':j\in{\bf Z}\}$ will be $\delta$-close to $\{f^j_{a,b}p:j\in{\bf
Z}\}$ if $(a',b')$ is close to $(a,b)$, since these sets are finite.   By Theorem
2, it follows that $S=c_{\cal G}({\rm Per}_N(a,b))$ for all $(a,b)$ in the homotopy.
Finally, since $S$ is finite, we conclude that the coding of each individual
periodic point is constant as we move $(a,b)$ through any point of the homotopy. 
In particular, the codings do not change as we move around $\gamma$, so $\rho_{\cal
G}(\gamma)$ is the identity transformation.

\bigskip\noindent{\bf \S3.  Proof of Theorem 2.}   We will show that an explicit automorphism of the shift, constructed by E. Brown, is not in $\Gamma$.  Brown uses the following scheme to describe automorphisms of the 2-shift.  Let us write ${\bf Z}_2=\{{\bf 0}, {\bf 1}\}$ for the field of two elements.  For
$x\in{\bf Z}_2$, let
$\bar x$ denote the element different from $x$.  Consider the polynomial $F:{\bf
Z}^4_2\to {\bf Z}_2$ defined by $F(x_1,x_2,x_3,x_4)=x_3 + x_1\bar x_2x_4$.  This
polynomial defines a map of $\Sigma_2$, which we will again denote as $F$, and which is given by $F:(x_n)\mapsto(y_n)$,
where $y_n=F(x_n,x_{n+1},x_{n+2}, x_{n+3})$.  Brown [B], shows that $F$ is an automorphism of $\Sigma_2$.  We let $F$ act on
$\Sigma_{\cal E}$.  Thus, if we apply
the formula for $F$, we have  $F({\bf 0}{\bf 0}{\bf 0}{\bf 1})^\infty= ({\bf 0}{\bf 1}{\bf 0}{\bf 0})^\infty$, and  $F({\bf 0}{\bf 0}{\bf 1}{\bf 1})^\infty= ({\bf 1}{\bf 1}{\bf 0}{\bf 1})^\infty$.

We will show that $F$ does not belong to the subgroup of $Aut(\Sigma_2)$ generated by
$\rho(\pi_1(HOV))$, $\rho(\pi_1({\cal W}\cap{\cal H^{\bf C}}))$ and the
shift.  If $(a,b)\in HOV\cap{\cal W}$, then each point of $J_{a,b}$ will have both an ${\cal E}$-coding and a ${\cal G}$-coding.  We will compare these two codings, using the graph $\hat{\cal G}$ in Figure 8.  

Let $P_4$ denote the set of twelve points of $\Sigma_2=\Sigma_{\cal E}$ with
period exactly four.  Points of period four do not belong to $W^s(P_\alpha)$, so they have unique ${\cal G}$-codings.  Thus we may represent them in the space $\Sigma_{\hat{\cal G}}$ as  $\pi_{\cal G}^{-1}P_4=X\cup Y$, where
$$ X = \left\{  {{\bf 0}\atop 0}{{\bf 0}\atop 0}{{\bf 0}\atop 1}{{\bf 1}\atop 2}\ ,\ {{\bf 1}\atop 2}{{\bf 1}\atop 1}{{\bf 1}\atop 2}{{\bf 0}\atop 1}\ ,\ {{\bf 0}\atop 0}{{\bf 0}\atop 0}{{\bf 1}\atop 1}{{\bf 1}\atop 2}\ ,\
 {{\bf 1}\atop 1}{{\bf 1}\atop 2}{{\bf 0}\atop 0}{{\bf 0}\atop 0}\ ,\ {{\bf 1}\atop 2}{{\bf 0}\atop 1}{{\bf 1}\atop 2}{{\bf 1}\atop 1}\ ,\
{{\bf 0}\atop 1}{{\bf 1}\atop 2}{{\bf 0}\atop 0}{{\bf 0}\atop 0}      \right\}$$ 
and
$$Y = \left\{ {{\bf 0}\atop 0}{{\bf 0}\atop 1}{{\bf 1}\atop 2}{{\bf 0}\atop 0}\ ,\ {{\bf 1}\atop 1}{{\bf 1}\atop 2}{{\bf 0}\atop 1}{{\bf 1}\atop 2}\ ,\ {{\bf 0}\atop 0}{{\bf 1}\atop 1}{{\bf 1}\atop 2}{{\bf 0}\atop 0}\ ,\
 {{\bf 1}\atop 2}{{\bf 0}\atop 0}{{\bf 0}\atop 0}{{\bf 1}\atop 1}\ ,\ {{\bf 0}\atop 1}{{\bf 1}\atop 2}{{\bf 1}\atop 1}{{\bf 1}\atop 2}\ ,\
{{\bf 1}\atop 2}{{\bf 0}\atop 0}{{\bf 0}\atop 0}{{\bf 0}\atop 1}\right\}.$$  
Of course ${{\bf 0}\atop 0}{{\bf 0}\atop 0}{{\bf 0}\atop 1}{{\bf 1}\atop 2}$ is just an admissible block of four symbols
from $\hat{\cal G}$; we use this as a notation to indicate the corresponding
periodic point which is obtained as the bi-infinite sequence in which
this block is repeated infinitely often.  And
$\pi_{\cal E}$ maps this periodic point to the element $({\bf 0}{\bf 0}{\bf 0}{\bf 1})^\infty\in\Sigma_{\cal E}$.

Now we consider the elements of the various subgroups and show that they respect the
sets $X$ and $Y$.  If 
$\sigma$ denotes the shift, then we have
$\sigma X=Y$ and $\sigma Y=X$.  The only nontrivial element of
$\rho\pi_1(HOV)$ is the automorphism $C\in Aut(\Sigma_{\cal E})$
which interchanges the symbols
${\bf 0}$ and ${\bf 1}$.  An inspection shows that $C(X)=X$ and $C(Y)=Y$.  

Now let $g$ be any
element of $\rho\pi_1({\cal W}\cap{\cal H^{\bf C}})$.  We see that $\pi_{\cal G}X=\{0012,
1200, 2121\}$ and $\pi_{\cal G}Y=\{0120,1212, 2001\}$.  By  Theorem 5, then, we must have
$gX=X$ and $gY=Y$.

It follows that if $g$ is in the span of these subgroups, it must preserve the set
of points of period  four, and we must have either
$gX=X$ or $gX=Y$.  On the other hand, we see from the calculation above that
$FX\cap X\ne\emptyset$ and
$FX\cap Y\ne\emptyset$, so $F$ cannot be in the span of these subgroups.

\bigskip\noindent{\bf \S4.  Type 3 maps and involutions.}  In this section we will establish Theorem 7, which leads to the proof of Theorem 1.   For a point of $HOV$ there is a Markov partition consisting of two pieces which allows us to construct a conjugacy $c$ from $K$ to $\Sigma_2$.  Note that $HOV\cap\{a<0,b\in{\bf R}\}$ is the subset of $HOV$ for which the (complex) horseshoe is contained in ${\bf R}^2$.  For a point $*\in HOV\cap\{a<0,b\in{\bf R}\}$,  one of the two Markov partition pieces lies to the left of the other, and we can assign to this  the symbol 0 and assign the symbol 1 to the other piece.  We denote this conjugacy by $c_*:K_*\to\Sigma_2$.   If we start at a point $*$ and move along a loop in $HOV$ which encircles the $b$-axis, then the two pieces of the partition are interchanged.  On the other hand, if we start at  $*$ and move around a loop inside $HOV$ which encircles the $a$-axis, then the pieces are not changed.  Thus we can take $*\in HOV\cap\{a<0,b\in{\bf R}\}$ as a canonical basepoint for $\pi_1({\cal H}^{\bf C})$.

Let us review the construction of $\rho:\pi_1({\cal H}^{\bf C},*)\to Aut(\Sigma_2)$.
Given a loop $\sigma$ based at $*$, we construct an automorphism $\rho(\sigma)$ of $\Sigma_2$ as follows. Structural stability gives a unique continuous family of conjugacies $h_t:K_{\sigma(t)}\to\Sigma_2$  such that $h_0=c_*$. We set $\rho(\sigma)=h_1\circ h_0^{-1}$. A loop which represents a generator of $\pi_1(HOV,*)$ is mapped by $\rho$ to the automorphism that interchanges the two symbols.

Let $\phi:(x,y)\mapsto(\bar x,\bar y)$ denote complex conjugation on the dynamical space ${\bf C}^2$.    Complex conjugation in parameter space commutes with $\phi$ in the sense:
$$\phi\circ f_{a,b}=f_{\bar a,\bar b}\circ\phi.$$
Since complex conjugation conjugates the dynamics, it preserves the set of points with bounded orbits, and so we may define the restriction $\phi_{a,b}:K_{a,b}\to K_{\bar a, \bar b}$.
Since $f_{a,b}$ and $f_{\bar a,\bar b}$ are smoothly conjugate the points $(a,b)$ is in ${\cal H}$ if and only if $(\bar a,\bar b)$ is in ${\cal H}$. Thus the sets ${\cal H}^{\bf C}$ and ${\cal H}^{\bf C}_0$ is invariant under complex conjugation applied to parameter space.

Now for a point  $(a,b)\in {\bf R}^2\cap {\cal H}^{\bf C}_0$ we let $\gamma:[0,1]\to{\cal H}^{\bf C}$ denote a path which starts at $*\in{\cal H}\cap HOV$ and ends at $(a,b)$. Let $\gamma^{-1}$ be the path obtained by parametrizing the path in the opposite direction as in the definition of the fundamental group.
Let $\bar\gamma$ be the path obtained by applying complex conjugation to the coordinates of $\gamma$.
Let $\hat\gamma=\gamma\cdot\bar\gamma^{-1}$ where the symbol $\cdot$ denotes concatenation of paths. We will investigate the correspondence between the point $(a,b)$ and the automorphism $\rho(\hat\gamma)$. The map $\rho$ is defined on loops. There is an analogous construction for paths which we now describe. If $\sigma:[0,1]\to {\cal H}^{\bf C}_0$ is a path then we there is a unique continuous path of topological conjugacies $g_t:K_{\sigma(0)}\to K_{\sigma(t)}$ such that $g_0$ is the identity. We then define $\alpha(\sigma)=g_1:K_{\sigma(0)}\to K_{\sigma(1)}$. The correspondence $\alpha$ satisfies properties analogous to those of $\rho$. In particular $\alpha(\sigma\cdot\tau)=\alpha(\tau)\circ\alpha(\sigma)$ and $\alpha(\sigma^{-1})=\alpha^{-1}(\sigma)$. If paths $\sigma$ and $\sigma'$ have the same endpoints and are homotopic relative to the endpoints then $\alpha(\sigma)=\alpha(\sigma')$.
When $\sigma$ is a loop we can describe $\rho(\sigma)$ in terms of $\alpha(\sigma)$ as follows: $\rho(\sigma)=c_*\circ \alpha(\sigma)\circ c_*^{-1}$.

The correspondence $\alpha$ is well behaved with respect to complex conjugation. We have
$$\phi_{\sigma(1)}\circ\alpha(\sigma)=\alpha(\bar\sigma)\circ\phi_{\sigma(0)}.$$
We can prove this by observing that if $g_t:K_{\bar\sigma(0)}\to K_{\bar\sigma(t)}$ is a family of conjugacies such that $g_0$ is the identity then $h_t=\phi^{-1}_{\sigma(t)}\circ g_t\circ\phi_{\sigma(0)}:
K_{\sigma(0)}\to K_{\sigma(t)}$ is a family of conjugacies with the property that $h_t$ is the identity. Thus $\alpha(\sigma)=h_1$ and $\alpha(\bar\sigma)=g_1$. So by setting $t=1$ in the above equation we have: $\alpha(\sigma)=\phi^{-1}_{\sigma(t)}\circ \alpha(\bar\sigma)\circ\phi_{\sigma(0)}$.

We are interested in the loop $\hat\gamma=\gamma\cdot\bar\gamma^{-1}$ based at $*$. It is convenient to introduce a closely related loop $\check\gamma=\bar\gamma^{-1}\cdot\gamma$. We can think of $\check\gamma$  as an alternate parametrization of the circle described $\hat\gamma$ which starts at the point $(a,b)$ rather than at $*$.
In terms of compositions of paths we have $\hat\gamma=\gamma\cdot\check\gamma\cdot\gamma^{-1}$.

Now applying the complex conjugation formula  to $\gamma$ we get $\phi_{a,b}\circ\alpha(\gamma)=\alpha(\bar\gamma)\circ\phi_*$. We note that since $K_*\subset{\bf R}^2$, the map $\phi_*$ is the identity. Thus we get
$$\alpha(\gamma)\circ\alpha(\bar\gamma^{-1})=\phi_{a,b}.$$
This gives us $\alpha(\check \gamma)=\phi_{a,b}.$

\proclaim Lemma 6.  If $(a,b)\in{\cal H}$, then $\rho(\hat\gamma)$ is the identity automorphism.  Otherwise, if $(a,b)\in{\cal H}^{\bf C}_0-{\cal H}$, then $\rho(\hat\gamma)$ is an involution.  If $\kappa\in\pi_1({\cal H}^{\bf C},*)$, then $\rho(\bar\kappa)=\rho(\kappa)$.

\noindent{\bf Proof.}    Since $\alpha(\check\gamma)$ coincides with complex conjugation, it is the identity transformation if $K_{a,b}\subset{\bf R}^2$, and  an element of order 2 otherwise.  Since $\alpha(\hat\gamma)=\alpha(\gamma\cdot\check\gamma\cdot\gamma^{-1})=\alpha(\gamma)^{-1}\circ\alpha(\check \gamma)\circ\alpha(\gamma)$ it follows that $\alpha(\hat\gamma)$ is topologically conjugate to $\alpha(\check\gamma)$ and that $\rho(\hat\gamma)$ is topologically conjugate to $\phi_{a,b}$. 

The condition that $(a,b)$ is a real horseshoe parameter is equivalent to $K_{a,b}\subset{\bf R}^2$.  Now $\rho(\hat\gamma)$, being conjugate to $\phi_{a,b}$, is also of order 1 or 2.

Let $\kappa$ be a loop based at $*$. Then $\phi_*\circ\alpha(\kappa)=\alpha(\bar\kappa)\circ\phi_*$. By the above observation $\phi_*$ is the identity so $\alpha(\kappa)=\alpha(\bar\kappa)$. It follows that $\rho(\kappa)=\rho(\bar\kappa)$.

 \medskip
Let us fix a point $(a,b)\in{\bf R}^2\cap{\cal H}^{\bf C}_0-{\cal H}$ and a path $\mu_0$ running from $*$ to  $(a,b)$.  All such paths are of the form  $\mu=\nu\cdot\mu_0$, where $\nu\in\pi_1({\cal H}^{\bf C}_0,*)$ is a loop with base point $*$.   It follows from Lemma 6 that the set of $\rho(\hat\mu)=\rho({\bar\nu}^{-1})\rho(\mu_0)\rho(\nu)$ for all paths $\mu$ is the conjugacy class of involutions 
$$C_{a,b}=\{g^{-1}\rho(\hat\mu_0)g:g\in\Gamma\}\subset Aut(\Sigma_2).$$

Now let us consider a pair of points $(a_j,b_j)\in{\bf R}^2\cap{\cal H}^{\bf C}_0$, $j=1,2$, and let $\gamma_j\subset{\cal H}^{\bf C}_0$ be arcs running from $*$ to  $(a_j,b_j)$.  We consider the conjugacy
$$\alpha:=\alpha(\gamma_{2})\alpha(\gamma_{1})^{-1}:K_{a_{1},b_{1}}\to K_{a_{2},b_{2}}.$$

\proclaim Theorem 7.  If $\rho(\hat\gamma_1)=\rho(\hat\gamma_2)$, then $\phi\circ\alpha=\alpha\circ\phi$.  It follows that $\alpha$ induces a conjugacy between $(f_{a_1,b_1},K^{\bf R}_{a_1,b_1})$ and $(f_{a_2,b_2},K^{\bf R}_{a_2,b_2})$.

\noindent{\bf Proof. }  By the discussion above, we have $\rho(\hat\gamma_{j})=\alpha(\gamma_{j})^{-1}\phi\alpha(\gamma_{j})$.  Setting these expressions equal, $\rho(\hat\gamma_{1})=\rho(\hat\gamma_{2})$, we find $\alpha(\gamma_{2})\alpha(\gamma_{1})^{-1}\circ\phi=\phi\circ\alpha(\gamma_{2})\alpha(\gamma_{1})^{-1}$.  
\bigskip
\noindent{\bf Proof of Theorem 1. }  Let $(a_{1},b_{1})$ and $(a_{2},b_{2})$ be two points of ${\bf R}^{2}\cap{\cal H}^{\bf C}_{0}$ such that $(f_{a_{1},b_{1}},K_{a_{1},b_{1}}^{\bf R})$ is not conjugate to  $(f_{a_{2},b_{2}},K_{a_{2},b_{2}}^{\bf R})$.  By Theorem 7, it follows that $\rho(\hat\gamma_{1})\ne\rho(\hat\gamma_{2})$.  Thus if $(a_{j},b_{j})$, $j=1,\dots,N$ correspond to distinct real conjugacy classes, then $\rho(\hat\gamma_{j})$, $j=1,\dots,N$ are distinct elements of $\Gamma$.

\bigskip
\noindent{\bf Remark.}  Figures 1 and 2 give evidence for the existence of nontrivial loops inside the slice of the horseshoe locus ${\cal H}^{\bf C}\cap\{b=.2\}.$  It is natural to ask whether such loops would also be nontrivial inside $\pi_{1}({\cal H}^{\bf C})$.  To address this, we consider an arc $\tau$ inside the complex slice of ${\cal H}^{\bf C}_{0}$, which starts and ends on the real $a$-axis.  It follows that $\hat\tau$ is contained in the slice, and  if $\tau$ connects two nonequivalent points of type 3, then  by Theorem 7, $\hat\tau$ generates a nontrivial element of $\pi_{1}({\cal H}^{\bf C})$.

\bigskip\noindent{\bf \S5.  Additional facts about $\rho$, and the shift locus.}  
There is another interesting homomorphism from $\pi_1({\cal H^{\bf C}},(a_0,b_0))$. We obtain this by considering the function $\Lambda$ on parameter space, where $\Lambda(a,b)$ is the larger Lyapunov exponent of $f_{a,b}$ with respect to the equilibrium measure. In [BS3] we showed that $\Lambda$ is a pluri-subharmonic function on parameter space which is pluri-harmonic in any region where the topological conjugacy type of $f_{a,b}$ is locally constant.  In terms of the operators $d$ and $d^c$, the condition that $\Lambda$ is pluri-harmonic is equivalent to the condition that $dd^c\Lambda=d(d^c\Lambda)=0$.  In particular this holds in the region ${\cal H^{\bf C}}$,  so $d^c\Lambda$ is a closed form that determines a real valued cohomology class $\theta$.

Consider the loop in the parameter region $HOV$ given by $\gamma(t)=(Re^{2\pi i t}, b)$ where $b\ne0$ and $R$ is sufficiently large. Let $[\gamma]$ be the corresponding homology class in ${\cal H^{\bf C}}$.

\proclaim Propositon 8. The cohomology class $\theta$ is non-zero. In particular $\theta([\gamma])=2\pi$.

\noindent Proof. A priori $\Lambda$ is defined on the set of parameters where $b\ne0$. It is shown in [BS3] that $\Lambda$ extends to a subharmonic function on all of ${\bf C}^2$ and that $\Lambda(a,0)=h(a)$ is given by the  function $h(a)$ which measures the rate of escape the critical point of $f(z)=z^2+a$. It follows that this extension is harmonic in the region where $|a|>2$. Thus we can evaluate $\theta([\gamma])$ by evaluating $d^c h$ on the loop $t\mapsto Re^{2\pi i t}$. Since $h$ approaches $\log |a|$ as $|a|\to\infty$ we get the value $2\pi$.

\proclaim Corollary 9. The homomorphism $\rho:\pi_1({\cal H}_0^{\bf C})\to Aut(\Sigma_2)$ is not injective.

\noindent Proof. We see that the loop represented by $\gamma$ has infinite order in $\pi_1({\cal H}_0^{\bf C})$ but maps to an element of order 2.

In connection with horseshoes, one may also consider the mappings which are conjugate to the shift.  Let us define the shift locus ${\cal S}$ in parameter space by the condition that $(a,b)\in{\cal S}$ if $J_{a,b}=K_{a,b}$, and the restriction of $f$ to this set is conjugate to the 2-shift $\Sigma_2$.  It follows that ${\cal H}^{\bf C}\subset{\cal S}$, but we do not know whether this containment is strict, since we do not know whether a map in the shift locus is necessarily hyperbolic. 

\proclaim Proposition 10.  ${\cal H}={\bf R}^2\cap{\cal H}^{\bf C}={\bf R}^2\cap{\cal S}$.

\noindent{\bf Proof.}  The only thing that needs to be proved is that if $(a,b)\in{\bf R}^2\cap{\cal S}$, then  $f_{a,b}$ is hyperbolic.  Since $f$ is quadratic and conjugate to the 2-shift, it is a real mapping of maximal entropy.  Since it is also topologically expansive, it follows from [BS8] that  $f$ is also hyperbolic.
\medskip
It is not known whether ${\cal S}$ is connected.  We let ${\cal S}_0$ denote the component of  ${\cal S}$ which contains $HOV$.   We note that any path $\gamma:[0,1]\to{\cal S}$ and any conjugacy $h_0:J_{\gamma(0)}\to\Sigma_2$, there is a unique, continuous family of conjugacies $h_t:J_{\gamma(t)}\to\Sigma_2$.  Thus we have a representation $\rho:\pi_1({\cal S}_0)\to Aut(\Sigma_2)$.  Further, it follows from  [BLS] that $\Lambda$ is pluriharmonic on ${\cal S}$, so the 1-form $\theta$ defined above defines a cohomology class on ${\cal S}$.
Thus, analogues of the results and conjectures in the previous sections apply to ${\cal S}$ as well as ${\cal H}^{\bf C}$.   

\proclaim Proposition 11.  Let complex parameters $(a,b)\in{\cal P}^2$ be given, and suppose that there is an invariant set $X\subset{\bf C}^2$ such that $(f_{a,b},X)$ is conjugate to $(\sigma,\Sigma_2)$.  Then $(a,b)\in{\cal S}$.

\noindent{\bf Proof.}  Let $\varphi:\Sigma_2\to X$ be a conjugacy, and let $\nu$ be the measure of maximal entropy for $\sigma$.  It follows from [BLS] that $\varphi_*\nu=\mu$.  Thus $X=J^*$.  Since this set is totally disconnected, it follows from Corollary 2.4 of [BS3] that $X=J=K$.

\bigskip
\bigbreak

\centerline{\bf References}
\medskip

\item{[BLS]} E. Bedford, M. Lyubich and J. Smillie,  Distribution of periodic points of polynomial diffeomorphisms of ${\bf C}^2$.  Invent. Math. 114 (1993), no. 2, 277--288.
\item{[BS3]} E. Bedford and J. Smillie, Polynomial diffeomorphisms of ${\bf C}^2$. III. Ergodicity, exponents and entropy of the equilibrium measure.  Math. Ann. 294 (1992), no. 3, 395--420.
\item{[BS8]} E. Bedford and J. Smillie,  Polynomial diffeomorphisms of
${\bf C}^2$.  VIII: Quasi-expansion.  American J. of Math., 124, 221--271,
(2002).
\item{[BS{\it i}]}  E. Bedford and J. Smillie,  Real polynomial diffeomorphisms with
maximal entropy: Tangencies.  Annals of Math., 160, 1--26 (2004).  
{\tt
www.arXiv.org/math.DS/0103038 }
\item{[BS{\it ii}]}  E. Bedford and J. Smillie,  Real polynomial
diffeomorphisms with maximal entropy: II. Small Jacobian, 
{\tt
www.arXiv.org/math.DS/0402355  }
\item{[BS]} E. Bedford and J. Smillie, A symbolic characterization of the horseshoe locus in the H\'enon family, in preparation.
\item{[BDK1]}  P. Blanchard, R. Devaney, and L. Keen, The dynamics of complex
polynomials and automorphisms of the shift. Invent. Math. 104 (1991), no. 3,
545--580. 
\item{[BDK2]}  P. Blanchard, R. Devaney, and L. Keen, Complex dynamics and symbolic dynamics, Proc. Symp. in Appl. Math.
\item{[B]} E. Brown, Periodic seeded arrays and automorphisms of the shift.
Trans. Amer. Math. Soc. 339 (1993), no. 1, 141--161.
\item{[DMS]} M.J. Davis, R. S. MacKay, and A. Sannami,  Markov shifts in the
H\'enon family.  Phys. D 52 (1991), no. 2-3, 171--178.
\item{[DN]} R. Devaney, Z. Nitecki,
Shift automorphisms in the H\'enon mapping. 
Comm. Math. Phys. 67 (1979), no. 2, 137--146.
\item{[EM]} H. El Hamouly and C. Mira, Lien entre les propri\'et\'es d'un
endomorphisme de dimension un et celles d'un diff\'eomorphisme de dimension
deux. (French) [Relation between the properties of a one-dimensional
endomorphism and those of a two-dimensional diffeomorphism] C. R. Acad. Sci.
Paris S\'er. I Math. 293 (1981), no. 10, 525--528.
\item{[HW]} P. Holmes and D. Whitley,  Bifurcations of one- and
two-dimensional maps. Philos. Trans. Roy. Soc. London Ser. A 311 (1984), no.
1515, 43--102. 
\item{[H]} ÊJ.H. Hubbard, The H\'enon mapping in the complex domain. Chaotic
dynamics and fractals (Atlanta, Ga., 1985), 101--111, Notes Rep. Math. Sci. Engrg.,
2, Academic Press, Orlando, FL, 1986.
\item{[HO]} J.H. Hubbard and R. Oberste-Vorth,  H\'enon mappings in the
complex domain II: Projective and inductive limits of polynomials, in: {\sl
Real and Complex Dynamical Systems}, B. Branner and P. Hjorth, eds. 89--132
(1995).
\item{[HP]}  J.H. Hubbard and K. Papadantonakis, Exploring the parameter
space of complex H\'enon mappings, to appear 
\item{[K]} B. Kitchens, {\sl Symbolic Dynamics}, Springer Verlag, 1998.
\item{[MNTU]} S. Morosawa, Y. Nishimura, M. Taniguchi, and T. Ueda, {\sl
Holomorphic Dynamics}, Cambridge U. Press, 2000.
\item{[Ob]} R. Oberste-Vorth, Complex horseshoes and the dynamics of 
mappings of two
complex variables, Ph.D dissertation, Cornell Univ., Ithaca, NY, 1987.
\item{[Ol]} R. Oliva, On the combinatorics of external rays in the dynamics
of the complex H\'enon map, Thesis, Cornell University, 1997. 
{\tt
www.math.cornell.edu/$\sim$dynamics}

\bigskip
\rightline{Indiana University}

\rightline{Bloomington, IN 47405}
\bigskip
\rightline{Cornell University}

\rightline{Ithaca, NY 14853}

\bye